\documentclass[a4paper,10pt,leqno]{article}
\usepackage{amsmath, amsfonts, amssymb, amsthm}
\numberwithin{equation}{section}

\begin{document}

\title{Global solution to liquid crystal flows in three dimensions}

\author{Wenke Tan\footnote{email:
tanwenkeybfq@163.com }\\ Department of Mathematics, Sun Yat-sen
University,\\510275 Guangzhou, China \smallskip\\ Zhaoyang Yin\footnote{email: mcsyzy@mail.sysu.edu.cn}\\
Department of Mathematics, Sun Yat-sen University,\\510275
Guangzhou, China }
\date{}
\maketitle

\begin{abstract}
In this paper, we mainly study a hydrodynamic system modeling the flow of nematic liquid crystals.
In three dimensions, we first establish local well-posedness of the initial-boundary value problem of the system. Then, we prove the existence of global strong solution to the system with small
initial-boundary condition.\\

\noindent \textit{Keywords}: Liquid crystal flow, local well-posedness, global strong solution.
\end{abstract}

\section{Introduction}
\par
In this paper we consider the following hydrodynamic system modeling the flow of liquid crystal materials
in three dimensions (see \cite{D,E,L-L,L-L-W})
\begin{equation}
 \left\{\begin{array}{ll}
\frac{\partial u}{\partial t}+u\cdot\nabla u-\nu\triangle u+\nabla P=-\lambda\nabla\cdot(\nabla d\odot\nabla d),
\ \ &in\ \  \mathbb{R}^{+}\times \Omega,\\
\frac{\partial d}{\partial t}+u\cdot\nabla d=\gamma(\triangle d+\mid\nabla d\mid^{2}d), \ \ &in\ \ \mathbb{R}^{+}\times \Omega,\\
\nabla\cdot u=0,\ \ \ \ \ \ \ \ \ \ \ \ \ \ \ \ \ \ \ \ \ \ \ &in\ \ \mathbb{R}^{+}\times \Omega, \\
\end{array}\right.
\end{equation} with initial-boundary conditions:
\begin{align}
(u(0,x),d(0,x))&=(u_0(x),d_0(x))\ \ \ x\in\Omega,\\
(u(t,x),d(t,x))&=(0,d_0(x))\ \ (t,x)\in\mathbb{R}^{+}\times\Omega.
\end{align}
Suppose that $\Omega\subseteq \mathbb{R}^3$ is a bounded smooth domain, $u(t,x):\mathbb{R}^{+}\times\Omega\rightarrow \mathbb{R}^3$ stands for the velocity field of the flow, $d(t,x):\mathbb{R}^{+}\times\Omega\rightarrow S^{2}$, the unit sphere in $\mathbb{R}^3$, is a unit-vector field that represents the macroscopic molecular orientation of the liquid crystal material and $P(t,x): \mathbb{R}^{+}\times \Omega\rightarrow \mathbb{R}$ is the pressure function.
The constants $\nu,\lambda,\gamma$ are positive constants that stand for the viscosity, the competition between kinetic energy and potential energy, and microscopic elastic relaxation time for the molecular orientation field. $\nabla d\odot\nabla d$ denotes the $3\times 3$ matrix whose the $(i,j)$ entry is given by $\nabla_i d\cdot\nabla_j d$ for $1\leq i,j\leq 3$. It is easy to see that $\nabla d\odot\nabla d=(\nabla d)^T\nabla d$, where $(\nabla d)^T$ denotes the transpose of matrix $\nabla d$.

System (1.1) is a simplified version of the Ericksen-Leslie model. General Ericksen-Leslie model reduces the Ossen-Frank model in the static case, for the hydrodynamics of nematic liquid crystals developed during the period from 1958 to 1968 \cite{D,E,L}. Since the general Ericksen-Leslie system is very complicated, we only study a simplified model of the Ericksen-Leslie system which can derive without destroying the basic structure. It is a macroscopic continuum description of the time evolution of the materials under the influence of both the flow field $u(t,x)$, and the macroscopic description of the microscopic orientation configurations $d(t,x)$ of rod-like liquid crystals. The system (1.1)-(1.3) is a system of the Navier-Stokes equation coupled with the harmonic map flows.

 In a series of papers, Lin \cite{Lin} and Lin-Liu \cite{L-L,L-L1} initiated the mathematical analysis of the system (1.1)-(1.3). Since the Erichsen-Leslie system (1.1)-(1.3) with $\mid d\mid=1$ is complicated, Lin and Liu \cite{L-L,L-L1} proposed to consider an approximation model of Ericksen-Leslie system by Ginzburg-Landau functional. More precisely, they replaced the Dirichlet functionals
 \begin{align*}
 \frac{1}{2}\int_{\Omega}\mid\nabla d\mid^2dx
 \end{align*} for $d:\Omega\rightarrow S^{n-1}$ by the Ginzburg-Landau functionals
 \begin{align*}
 \int_{\Omega}(\frac{1}{2}\mid\nabla d\mid^2+\frac{(1-\mid d\mid^2)^2}{4\epsilon})dx
 \end{align*} for $d:\Omega\rightarrow \mathbb{R}^n (\epsilon>0)$. In \cite{L-L}, Lin and Liu proved the global existence of solutions in dimensions two or three. In \cite{L-L1}, Lin and Liu proved partial regularity of weak solutions in dimension three. Furthermore, Lin and Liu in \cite{L-L2} proved existence of solutions for the general Ericksen-Leslie system and also analyzed the limits of weak solutions as $\epsilon\rightarrow 0$. In \cite{H-W}, Hu and Wang give the existence of global strong solution and prove that all the weak solutions constructed in \cite{L-L} must be equal to the unique strong solution.

Recently, Lin, Lin and Wang \cite{L-L-W} studied the system (1.1)-(1.3) in two dimensions. They established the global existence and partial regularity of the global weak solution and performed the blow-up analysis at each singular time. Hong \cite{H} proved the global existence of the system (1.1)-(1.3) in two dimensions independently.

The aim of this paper is to establish the short-time solution for general initial-boundary condition and the global existence for small initial-boundary
conditions for the system (1.1)-(1.3).\\

\noindent\textbf{Notations}\textit{\ In this paper, we denote
$W^{m,q}(\Omega)$ the set of function in $L^q(\Omega)$ whose
derivatives up to order $m$ belong to $L^q(\Omega)$.  For $T>0$
and a function space $X$, we denote by $L^p(0,T;X)$ the set of
Bochner measurable $X-value$ time dependent functions $f$ such
that $t\rightarrow \parallel f\parallel_{X}$ belong to $L^p(0,T)$.
The function  $u\in V^{1,0}_2(Q_T)$ is defined as
\begin{align*}
u\in C([0,T],L^2(\Omega))\cap W^{1,0}_2(Q_T),\\
\sup_{[0,T]}\parallel u(t,\cdot)\parallel_{L^2(\Omega)}+\parallel \nabla u\parallel_{L^2(Q_T)}\leq \infty,
\end{align*} where $Q_T=(0,T)\times \Omega$.
The space $D_{A_q}^{1-\frac{1}{p},p}$ represents some fraction domain of Stokes operator in $L^q$ (see Sect.2.3 in \cite{Dan}). Roughly, the vector-fields of $D_{A_q}^{1-\frac{1}{p},p}$ are vectors which have $2-\frac{2}{p}$ derivatives in $L^q(\Omega)$, are divergence-free, and vanish on $\partial \Omega$.  $B_{q,p}^{s}(\Omega)$ represents the Besov space \cite{B-L} which can be regarded as the interpolation space between $L^q(\Omega)$ and $W^{s+\epsilon,q}(\Omega)$. From Proposition 2.5 in \cite{Dan}, we can get
\begin{align}
D_{A^q}^{1-\frac{1}{p},p}\hookrightarrow B_{q,p}^{2(1-\frac{1}{p})}(\Omega)\cap L^q(\Omega).
\end{align}
   Since the space variables are in $\Omega$, if there is no ambiguity, we write $L^q(\Omega)$, $W^{m,q}(\Omega)$, $B_{q,p}^{s}(\Omega)$ as $L^q$, $W^{m,q}$, $B_{q,p}^{s}$
respectively.
}
\\

\noindent\textbf{Definition 1.1.}\textit {\ For $T>0$ and $1<p,q<\infty$,
we denote by $E_{T}^{p,q}$ the set of triplets (u, d, P) such that\\
$$u\in C(0,T; D_{A^{q}}^{1-\frac{1}{p},p})\cap L^p(0,T;W^{2,q}(\Omega)\cap W^{1,q}_0(\Omega)),$$ $$\partial_t u\in L^p(0,T;L^q(\Omega)), \nabla\cdot u=0,$$\\
$d\in C(0,T;B_{q,p}^{2(1-\frac{1}{p})})\cap L^p(0,T;W^{2,q}(\Omega))$, $\partial_t d\in L^p(0,T;L^q(\Omega))$,\\
$P\in L^p(0,T;W^{1,q}(\Omega))$, $\int_{\Omega}P dx=0$.\\
The corresponding norm is denoted by $\parallel \cdot\parallel_{E_{T}^{p,q}}$.
\begin{align*}
\parallel (u,d,P)\parallel_{E_{T}^{p,q}}=&\sup_{[0,T]}\parallel u \parallel_{D_{A^q}^{1-\frac{1}{p},p}}+\sup_{[0,T]}\parallel d \parallel_{B_{q,p}^{2(1-\frac{1}{p})}}+\parallel u \parallel_{L^p(0,T;W^{2,q})}\\&+\parallel \partial_t u\parallel_{L^p(0,T;L^q)}+\parallel d \parallel_{L^p(0,T;W^{2,q})}+\parallel \partial_t d \parallel_{L^p(0,T;L^q)}.\notag
\end{align*}
}\\

\par
 Our main results can be stated as follows.
\\

\noindent\textbf{Theorem 1.1.}\textit{\ Let $\Omega$ be a smooth bounded domain in $\mathbb{R}^3$ and $(1-\frac{2}{p})\cdot q>3$. If $u_0\in D_{A^q}^{1-\frac{1}{p},p}$ and $d_0\in B_{q,p}^{2(1-\frac{1}{p})}\cap C^{2,\alpha}(\partial \Omega),$ then
\\
\\
(1) There exists a $T_0>0,$ such that the system (1.1) with the initial-boundary condition (1.2)-(1.3) has a unique local strong solution $(u,d,P)\in E_{T_0}^{p,q}$ in $(0,T_0)\times \Omega$. Moreover the solution continuously depends on initial data.
\\
\\
(2) For any given unit vector $e\in S^2$,there exists a $\delta>0$, such that, if the initial data satisfies $d_0|_{\partial \Omega}=e$ and
\begin{align*}
\parallel u_0 \parallel_{D_{A^q}^{1-\frac{1}{p},p}}+\parallel d_0-e \parallel_{B_{q,p}^{2(1-\frac{1}{p})}\cap C^{2,\alpha}(\partial \Omega)}\leq \delta,
\end{align*} then the system (1.1)-(1.3) has a unique global strong solution $(u,d,P)\in E_{T}^{p,q}$ in $(0,T)\times\Omega$ for all $T>0.$
}\\

We obtain our results in the spirit of \cite{Dan,L-L-W}. The main difficulty is the low integrability of nolinear item $\nabla\cdot(\nabla d\odot\nabla d)$. To overcome this problem, we add a certain condition to the initial data such that we can get the estimate of the $\parallel \nabla d\parallel_{L^\infty(0,T;L^\infty)}$.
The paper is written as follows. In Section 2, we give some useful lemmas.  In Section 3, we prove the local well-posedness. In Section 4, we prove the global existence.

\noindent\textbf{Remark 1.1}\textit{\ In this paper we only prove the results in three dimensions, but we point out that our method can deal with the system (1.1)-(1.3) in higher dimensions.}
\section{Preliminaries }
\par

In this section, we give some useful lemmas which will be used in the sequel.

\noindent\textbf{Lemma 2.1} \cite{A}\textit{\ Given $1<p,q<\infty,$ $u_0\in B_{q,p}^{2(1-\frac{1}{p})}$
and $f\in L^p(0,T;L^q)$}. Then the Cauchy problem
\begin{align}
\frac{\partial u}{\partial t}-\triangle u=f, \ \ u|_{t=0}=u_0,
\end{align}
has a unique solution $u$ satisfying
\begin{align}
&\parallel u\parallel_{W^{1,p}(0,T; L^q)}+\parallel u\parallel_{L^p(0,T;W^{2,q})}\\\leq &C_1(\parallel f\parallel_{L^p(0,T;L^q)}+\parallel u_0\parallel_{B_{q,p}^{2(1-\frac{1}{p})}}),\notag
\end{align} where $C_1$ is independent of $u_0,$ $f,$ and $T$. Moreover, there exists a positive constant $C_2$ independent of $f$ and $T$ such that
\begin{align}
\sup_{t\in(0,T)}\parallel u\parallel_{B_{q,p}^{2(1-\frac{1}{p})}}\leq C_2(\parallel f\parallel_{L^p(0,T;L^q)}+\parallel u_0\parallel_{B_{q,p}^{2(1-\frac{1}{p})}}).
\end{align}

\noindent\textbf{Lemma 2.2}  \cite{Dan}\textit{\ Let $\Omega$ be a $C^{2+\epsilon}$ bounded domain in $\mathbb{R}^N$ and $1<q,p<\infty$. Assume that $u_0\in D_{A^q}^{1-\frac{1}{p},p}$ and $f\in L^p(\mathbb{R}^{+},L^q).$ Then the system
\\\begin{equation}
\left\{\begin{array}{ll} \frac{\partial u}{\partial t}-\triangle u+\nabla P=f,\ \ \int_{\Omega}P dx=0,&\\
\nabla\cdot u=0,\ \ u|_{\partial\Omega}=0 ,\\
 u|_{t=0} = u_0,
\end{array}\right.
\end{equation}
has a unique solution $u,P$ satisfying the following inequality for all $T\geq 0$:
\begin{align}
&\parallel u(T)\parallel_{D_{A^q}^{1-\frac{1}{p},p}}+(\int_0^T\parallel (\nabla P,\nabla^2 u,\partial_t u)\parallel_{L^q}^p dt)^{\frac{1}{p}}\\\leq &C_3(\parallel u_0\parallel_{D_{A^q}^{1-\frac{1}{p},p}}+(\int_0^T\parallel f(t)\parallel_{L^q}^p)^{\frac{1}{p}}),\notag
\end{align}
with $C_3=C(q,p,N,\Omega)$.}
\\ Using Lemma 2.1, we can prove there is a similar conclusion for initial-boundary problem.
\\

\noindent\textbf{Theorem 2.1}\textit{\ Let $\Omega$ be a bounded smooth domain in $\mathbb{R}^3$ and $(1-\frac{2}{p})\cdot q>3$. If $u_0\in B_{q,p}^{2(1-\frac{1}{p})}\cap C^{2,\alpha}(\partial \Omega)$ and $f\in L^p(0,T;L^q)$, then the initial-boundary problem
\begin{align}
\frac{\partial u}{\partial t}-\triangle u=f,\ \ \ u|_{\partial {p}Q_T}=u_0
\end{align} has a unique solution $u$ satisfying
\begin{align}
&\parallel u\parallel_{L^\infty(0,T;B_{q,p}^{2(1-\frac{1}{p})})}+\parallel u\parallel _{L^p(0,T;W^{2,q})}+\parallel \partial_t u\parallel_{L^p(0,T;L^q)}\\\leq & C(\parallel f\parallel_{L^p(0,T;L^q)}+\parallel u_0\parallel_{B_{q,p}^{2(1-\frac{1}{p})}}+T^{\frac{1}{p}}\parallel u_0\parallel_{C^{2,\alpha}(\partial \Omega)}),\notag
\end{align} where $\partial_ {p}{Q_T}=(0,T)\times\partial \Omega\cup\{0\}\times \Omega$.
}
\begin{proof}
Since $u_0\in C^{2,\alpha}(\partial \Omega)$, by the standard elliptic theory we get that there exists a unique solution $u^1\in C^{2,\alpha}(\bar\Omega)$ satisfies
\begin{align}
\triangle u^1=0,\ \ \ u^1|_{\partial \Omega}=u_0.
\end{align} It is clear that $u_0-u_1\in {B_{q,p}^{2(1-\frac{1}{p})}}$ and $(u_0-u_1)|_{\partial \Omega}=0$. Using Lemma 2.1, the equation
\begin{align*}
\frac{\partial u}{\partial t}-\triangle u=f, \ \ u|_{t=0}=u_0-u_1,
\end{align*} has a solution $u_2\in B_{q,p}^{2(1-\frac{1}{p})}$ and $u_2(t,x)|_{\partial \Omega}=0$.
A direct computation shows that $u_1+u_2$ is a solution of the initial-boundary problem (2.6). Using lemma 2.1 and Schauder's estimate, we deduce (2.7).
\end{proof}

\noindent\textbf{Theorem 2.2}\textit{\ Let $(1-\frac{2}{p})\cdot
q>3$ and $u\in L^\infty(0,T;D_{A^q}^{1-\frac{1}{p},p})$. If $d\in
L^p(0,T;W^{2,q})\cap W^{1,p}(0,T;L^q)\cap
L^\infty(0,T;B_{q,p}^{2(1-\frac{1}{p})})$ is a solution of the
following nonlinear parabolic problem
\begin{align}
\frac{\partial d}{\partial t}-\triangle d-\mid\nabla d\mid^2d+u\cdot\nabla d&=0,\ \ in \ \ (0,T)\times \Omega,\\
d|_{\partial p_{Q_T}}&=d_0,
\end{align} where $d_0: \Omega\rightarrow S^2$. Then, $\mid d\mid=1$ in $[0,T)\times \Omega$.
}
\begin{proof}
Multiplying (2.9) by $d$, we get
\begin{align}
\frac{\partial (\mid d\mid^2-1)}{\partial t}-\triangle(\mid d\mid^2-1)-2\mid \nabla d\mid^2(\mid d\mid^2-1)-u\cdot\nabla(\mid d\mid^2-1)=0.
\end{align} By the assumption $(1-\frac{2}{p})\cdot q>3$ and (1.4), we have
\begin{align}
\parallel \nabla d\parallel_{L^\infty(0,T;L^\infty)} \leq\parallel \nabla d\parallel_{L^\infty(0,T;B_{q,p}^{1-\frac{2}{p}})}\leq\parallel d\parallel_{L^(0,T;B_{q,p}^{2(1-\frac{1}{p})})},\\
\parallel u\parallel_{L^\infty(0,T;L^\infty)}\leq \parallel u\parallel_{L^\infty(0,T; B_{q,p}^{2(1-\frac{1})})}\leq\parallel u\parallel_{L^\infty(0,T;D_{A^q}^{1-\frac{1}{p},p})}.
\end{align} Noticing that $(1-\frac{2}{p})\cdot q>3$, we get
\begin{align*}
\parallel \nabla (\mid d\mid^2-1)\parallel_{L^2(Q_T)}\leq &C \parallel d\parallel_{L^\infty(Q_T)}\parallel\nabla d\parallel_{L^2(Q_T)}\\\leq &C\parallel d\parallel_{L^\infty(0,T;B_{q,p}^{2(1-\frac{1}{p})})}\parallel d\parallel_{L^p(0,T;W^{2,q})}.\notag
\end{align*} Since $\partial_t d\in L^p(0,T;L^q)$ and $d\in L^\infty(0,T;L^\infty)$, we can get that $$(\mid d\mid^2-1)\in C([0,T],L^2).$$
 This yields $(\mid d\mid^2-1)\in V^{1,0}_2(Q_T)$.
Thus, we have that the function $\mid d\mid^2-1$ satisfies the following equation
\begin{align}
\frac{\partial f}{\partial t}-\triangle f-2\mid \nabla d\mid^2f-u\cdot\nabla f &=0,\\
f|\partial P_{Q_T}&=0,
\end{align} in $V_2^{1,0}(Q_T)$. The uniqueness implies $\mid d\mid^2-1=0$ in $Q_T$. This proves the lemma.
\end{proof}

\section{Local well-posedness}
\par
 In this section, we prove the local well-posedness of the system (1.1) with the initial boundary value (1.2)-(1.3).

Noticing that $\nabla P=\nabla (P-\int_{\Omega}Pdx)$, we can assume that
\begin{align*}
\int_{\Omega} Pdx=0.
\end{align*}
Since the exact values of $\nu,\lambda,\gamma$ do not play a role, we henceforth assume$$\nu=\lambda=\gamma=1.$$ Thus, we can rewrite the system (1.1) as
\begin{equation}
 \left\{\begin{array}{ll}
\frac{\partial u}{\partial t}+u\cdot\nabla u-\triangle u+\nabla P=-\nabla\cdot(\nabla d\bigodot\nabla d),
\ \ &in\ \  \mathbb{R}^{+}\times \Omega,\\
\frac{\partial d}{\partial t}+u\cdot\nabla d=(\triangle d+\mid\nabla d\mid^{2}d), \ \ &in\ \ \mathbb{R}^{+}\times \Omega,\\
\nabla\cdot u=0,\ \ \ \ \ \ \ \ \ \ \ \ \ \ \ \ \ \ \ \ \ \ \ &in\ \ \mathbb{R}^{+}\times \Omega, \\
\int_{\Omega} Pdx=0,
\end{array}\right.
\end{equation} with initial-boundary conditions
\begin{align}
\left\{\begin{array}{ll}
(u(0,x),d(0,x))=(u_0(x),d_0(x))\ \ x\in\Omega,\\
(u(t,x),d(t,x))=(0,d_0(x)),\ \ \ \ \ (t,x)\in\mathbb{R}^{+}\times\partial\Omega.
\end{array}\right.
\end{align}

Now, we prove the local existence.

Firstly, we linearize the system (3.1)-(3.2) and construct approximate solutions. Set $(u^0(t,x),d^0(t,x))=(u_0,d_0)$. Then given $(u^n,d^n,P^n)$ as the solution of\\

\begin{equation}
 \left\{\begin{array}{ll}
\frac{\partial u^n}{\partial t}-\triangle u^n+\nabla P^n=-u^{n-1}\cdot\nabla u^{n-1}-\nabla\cdot(\nabla d^{n-1}\bigodot\nabla d^{n-1}),
\ \ &in\ \  \mathbb{R}^{+}\times \Omega,\\
\frac{\partial d^n}{\partial t}-\triangle d^n=-u^{n-1}\cdot\nabla d^{n-1}+\mid\nabla d^{n-1}\mid^{2}d^{n-1}, \ \ &in\ \ \mathbb{R}^{+}\times \Omega,\\
\nabla\cdot u^n=0,\ \ \ \ \ \ \ \ \ \ \ \ \ \ \ \ \ \ \ \ \ \ \ &in\ \ \mathbb{R}^{+}\times \Omega, \\
\int_{\Omega} P^ndx=0,
\end{array}\right.
\end{equation} with initial-boundary data
\begin{align}
\left\{\begin{array}{ll}
(u^n(0,x),d^n(0,x))=(u_0(x),d_0(x))\ \ x\in\Omega,\\
(u^n(t,x),d^n(t,x))=(0,d_0(x))\ \ \ \ \ (t,x)\in\mathbb{R}^{+}\times\partial\Omega.
\end{array}\right.
\end{align}
From Theorem 2.1 and Lemma 2.2, we can obtain that the sequence $\{(u,^n,d^n,P^n)\}_{n\in \mathbb{N}}$
belong to $E_{T}^{p,q}$ for any $T>0$. Moreover we have the following estimate between $(u^n,d^n,P^n)$ and $(u^{n-1},d^{n-1},P^{n-1})$.
\begin{align}
&\parallel u^{n+1}\parallel_{L^\infty(0,T;D_{A^q}^{1-\frac{1}{p},p})}+(\int_{0}^{T}\parallel (\nabla P^{n+1},\nabla^2 u^{n+1},\partial_t u^{n+1})\parallel^p_{L^q}dt)^{\frac{1}{p}}\\\leq &C(\parallel u_0\parallel_{D_{A^q}^{1-\frac{1}{p},p}}+(\int_{0}^{T}\parallel u^n\cdot\nabla u^n+\nabla\cdot(\nabla d^n\odot\nabla d^n)\parallel^p_{L^q}dt)^{\frac{1}{p}}),\notag
\end{align}
\begin{align}
&\parallel d^{n+1}\parallel_{L^\infty(0,T;B_{q,p}^{2(1-\frac{1}{p})})}+\parallel d^{n+1}\parallel_{L^p(0,T;W^{2,q})}+\parallel \partial_t d^{n+1}\parallel_{L^{p}(0,T;L^q)}\\\leq &C(\parallel d_0\parallel_{B_{q,p}^{2(1-\frac{1}{p})}}+(\int_{0}^{T}\parallel -u^n\cdot\nabla d^n+\mid\nabla d^n\mid^2d^n\parallel^p_{L^q}dt)^{\frac{1}{p}}+T^{\frac{1}{p}}\parallel d_0\parallel_{C^{2,\alpha}(\partial \Omega)}).\notag
\end{align} Note that $u^{n+1}|_{\partial \Omega}=0$. Then the divergence theorem implies
\begin{align*}
\int_\Omega \nabla u^{n+1}dx=0.
\end{align*} Thus, using Poincar\'{e}'s inequality, we get
\begin{align*}
\parallel u^{n+1}\parallel_{W^{2,q}}\leq C \parallel \nabla^2 u^{n+1}\parallel_{L^q}.
\end{align*} Therefore (3.5) can be written as
\begin{align}
&\parallel u^{n+1}\parallel_{L^\infty(0,T;D_{A^q}^{1-\frac{1}{p},p})}+(\int_{0}^{T}\parallel (\nabla P^{n+1},u^{n+1},\nabla^2 u^{n+1},\partial_t u^{n+1})\parallel^p_{L^q}dt)^{\frac{1}{p}}\\\leq &C(\parallel u_0\parallel_{D_{A^q}^{1-\frac{1}{p},p}}+(\int_{0}^{T}\parallel u^n\cdot\nabla u^n+\nabla\cdot(\nabla d^n\odot\nabla d^n)\parallel^p_{L^q}dt)^{\frac{1}{p}}).\notag
\end{align}

Secondly, we present a uniform estimate for the sequence $\{(u^n,d^n,P^n)\}_{n\in \mathbb{N}}$.
Define
\begin{align}
 F_n(t)=&\parallel u^n\parallel_{L^\infty(0,t;B_{q,p}^{2(1-\frac{1}{p})})}+\parallel u^n\parallel_{L^p(0,t;W^{2,q})}\\&+\parallel \partial_t u^n\parallel_{L^p(0,t;L^q)}+\parallel \nabla P^n\parallel_{L^p(0,t;L^q)}\notag,\\
 E_n(t)=&\parallel d^n\parallel_{L^\infty(0,t;B_{q,p}^{2(1-\frac{1}{p})})}+\parallel d^n\parallel_{L^p(0,t;W^{2,q})}+\parallel \partial_t d^n\parallel_{L^p(0,t;L^q)},\\
 H_n(t)=&F_n(t)+E_n(t),\\
  F_0=&\parallel u_0\parallel_{D_{A^q}^{(1-\frac{1}{p}),p}},\ \ E_0=\parallel d_0\parallel_{B_{q,p}^{2(1-\frac{1}{p})}},\ \ H_0=F_0+E_0.
\end{align}

\noindent\textbf{Lemma 3.1}\textit{\ Let $(1-\frac{2}{p})\cdot q>3$ and  $u_0(x)\in D_{A^q}^{1-\frac{1}{p},p},$ $d_0(x)\in B_{q,p}^{2(1-\frac{1}{p})}\cap C^{2,\alpha}(\partial \Omega)$. Then,
there exists a positive $T_0$ such that the sequence $\{(u^n,d^n,P^n)\}_{n\in\mathbb{N}}$ is uniformly bounded in $E_{T_0}^{p,q}$.
}
\begin{proof}
Using the fact
\begin{align}
\nabla \cdot(\nabla d\odot\nabla d)=&\nabla(\frac{\mid\nabla d\mid^2}{2})+\triangle d\cdot\nabla d\\=&\nabla^2d\cdot\nabla d+\triangle d\cdot \nabla d\notag
\end{align} and (2.12)-(2.13), we have that
\begin{align}
F_{n+1}(t)\leq &C(F_0+\parallel u^n\cdot\nabla u^n\parallel_{L^p(0,t;L^q)}\\&+\parallel\nabla\mid\nabla d^n\mid^2\parallel_{L^p(0,t;L^q)}+\parallel \triangle d^n\cdot\nabla d^n\parallel_{L^p(0,t;L^q)})\notag\\\leq &C(F_0+t^{\frac{1}{p}}\parallel u^n\parallel_{l^\infty(0,t;L^\infty)}\parallel \nabla u^n\parallel_{L^\infty(0,t;L^\infty)}\notag\\&+\parallel \nabla d\parallel_{L^\infty(0,t;L^\infty)}\parallel \triangle d^n\parallel_{L^{p}(0,t;L^q)})\notag\\\leq &C(F_0+t^{\frac{1}{p}}\parallel u^n\parallel^2_{L^\infty(0,t; D_{A^q}^{1-\frac{1}{p},p})}\notag\\&+ \parallel d^n\parallel_{L^\infty(0,t;B_{q,p}^{2(1-\frac{1}{p})})}\parallel d^n\parallel_{L^p(0,t;W^{2,q})})\notag\\\leq &C(F_0+t^{\frac{1}{p}}F^2_n(t)+E^2_n(t)),\notag
\end{align} due to (3.7). Similarly, we obtain from (3.6) that
\begin{align}
E_{n+1}(t)\leq& C(\parallel u^n\cdot \nabla d^n\parallel_{L^p(0,t;L^q)}+\parallel \mid\nabla d^n\mid^2 d^n\parallel_{L^p(0,t;L^q)})\\\leq &C(E_0+t^{\frac{1}{p}}\parallel u^n\parallel_{L^\infty(0,t;L^\infty)} \parallel \nabla d^n\parallel_{L^\infty(0,t;L^\infty)}\notag\\&+t^{\frac{1}{p}}\parallel \nabla d^n\parallel^2_{L^\infty(0,t;L^\infty)}\parallel d^n\parallel_{L^\infty(0,t;L^\infty)}+t^{\frac{1}{p}}\parallel d_0\parallel_{C^{2,\alpha}(\partial \Omega)})\notag\\\leq &C(E_0+t^\frac{1}{p}F_n(t)E_n(t)+t^{\frac{1}{p}}E^3_n(t)+t^{\frac{1}{p}}\parallel d_0\parallel_{C^{2,\alpha}(\partial \Omega)}).\notag
\end{align}

Combining (3.13) with (3.14) yields that
\begin{align}
H_{n+1}(t)\leq C\{H_0+t^\frac{1}{p}[F^2_n(t)+F_n(t)E_n(t)+E^3_n(t)+\parallel d_0\parallel_{C^{2,\alpha}(\partial \Omega)}]+E^2_n(t)\}.
\end{align} Plugging (3.14) into (3.15), we get that
\begin{align}
H_{n+1}(t)\leq &C\{H_0+t^\frac{1}{p}(F^2_n(t)+F_n(t)E_n(t)+E^3_n(t)+\parallel d_0\parallel_{C^{2,\alpha}(\partial \Omega)})\\
+&C^2(E_0+t^\frac{1}{p}F_{n-1}(t)E_{n-1}(t)+t^{\frac{1}{p}}E^3_{n-1}(t)+t^{\frac{1}{p}}\parallel d_0\parallel_{C^{2,\alpha}(\partial \Omega)})^2\}.\notag
\end{align} By a direct computation, we obtain
\begin{align}
H_{n+1}(t)\leq C&\{H_0+C^2H^2_0+(C^2+\parallel d_0\parallel_{C^{2,\alpha}(\partial \Omega)})t^\frac{2}{p}(H^6_{n-1}+H^5_{n-1}+H^4_{n-1})\\
+&t^\frac{1}{p}(H^3_n+H^2_n)+(2C^2H_0+\parallel d_0\parallel_{C^{2,\alpha}(\partial \Omega)})t^\frac{1}{p}(H^3_{n-1}+H^2_{n-1})\},\notag\\
+&t^{\frac{1}{p}}\parallel d_0\parallel_{C^{2,\alpha}(\partial \Omega)}.\notag
\end{align}  where the constant $C$ only depends on $\Omega$. Since $H_n(t)\rightarrow H_0$ as $t\rightarrow 0$, we can assume that there exists a $T>0$ such that
$H_n(t)\leq 2CK H_0$ and $H_{n-1}(t)\leq 2CK H_0$ on $[0,T]$ for some fixed $n$, where $K=1+C^2H_0$. Then,
\begin{align}
H_{n+1}(t)\leq &C\{(H_0+C^2H^2_0)+C^3t^\frac{2}{p}(2^6K^6H^6_0+2^5K^5H^5_0+2^4K^4H^4_0)\\
&+t^\frac{1}{p}(C+2C^2H_0+\parallel d_0\parallel_{C^{2,\alpha}(\partial \Omega)})(2^3K^3H^3_0+2^2K^2H^2_0)\notag\\&+\parallel d_0\parallel_{C^{2,\alpha}(\partial \Omega)}t^{\frac{2}{p}}(2^6K^6H^6_0+2^5K^5H^5_0+2^4K^4H^4_0)\},\notag\\=&CKH_0+C\{C^3t^\frac{2}{p}(2^6K^6H^6_0+2^5K^5H^5_0+2^4K^4H^4_0)\notag\\
&+t^\frac{1}{p}(C+2C^2H_0+\parallel d_0\parallel_{C^{2,\alpha}(\partial \Omega)})(2^3K^3H^3_0+2^2K^2H^2_0)+t^\frac{1}{p}\notag\\&+\parallel d_0\parallel_{C^{2,\alpha}(\partial \Omega)}t^{\frac{2}{p}}(2^6K^6H^6_0+2^5K^5H^5_0+2^4K^4H^4_0)\}\notag\\&+t^{\frac{1}{p}}\parallel d_0\parallel_{C^{2,\alpha}(\partial \Omega)}\notag
\end{align} If we choose $T_0\leq T$ such that
\begin{align}
&(C^3+\parallel d_0\parallel_{C^{2,\alpha}(\partial \Omega)})T_0^\frac{2}{p}(2^6K^6H^5_0+2^5K^5H^4_0+2^4K^4H^3_0)
\\&+T_0^\frac{1}{p}(C+2C^2H_0+\parallel d_0\parallel_{C^{2,\alpha}(\partial \Omega)})(2^3K^3H^2_0+2^2K^2H_0)\notag\\\leq &K-1,\notag
\end{align} and
\begin{align*}
T_0^{\frac{1}{p}}\parallel d_0\parallel_{C^{2,\alpha}(\partial \Omega)}\leq CH_0,
\end{align*}then
\begin{align}
H_{n+1}(T_0)\leq 2CK H_0.
\end{align} Since $H_{n+1}(t)$ is increasing, we get
\begin{align}
H_{n+1}(t)\leq 2CK H_0, \ \ \ on\ \ \ [0,T_0].
\end{align}
Arguing by induction, we deduce that $\{(u^n,d^n,P^n)\}_{n\in \mathbb{N}}$ is uniformly bounded in $E^{q,p}_{T_0}$.
\end{proof}
Next, we establish the convergence of this approximate solutions sequence $\{(u^n,d^n,P^n)\}_{n\in \mathbb{N}}$.

\noindent\textbf{Lemma 3.2}\textit{\ Let $(1-\frac{2}{p})\cdot q>3,$ $u_0\in D_{A_q}^{1-\frac{1}{p},p}$ and $d_0\in B_{q,p}^{2(1-\frac{1}{p})}\cap C^{2,\alpha}(\partial \Omega).$ Given solution sequence $\{(u^n,d^n,P^n)\}$ constructed in Lemma 3.1. Then, there
exists a positive constant $T_1\leq T_0$ such that $\{(u^n,d^n,P^n)\}_{n\in \mathbb{N}}$ converges in $E_{T_1}^{q,p}$.
}
\begin{proof}
Let
\begin{align}
D(u^n)\doteq u^{n+1}-u^n;\ \ D(d^n)\doteq d^{n+1}-d^n;\ \ D(P^n)\doteq P^{n+1}-P^n.
\end{align}
Define
\begin{align}
DF_n(t)=&\parallel D(u^n)\parallel_{L^\infty(0,t;D_{A^q}^{1-\frac{1}{p},p})}+\parallel D(u^n)\parallel_{L^p(0,t;W^{2,q})}\\&+\parallel\partial_t D(u^n)\parallel_{L^p(0,t;L^q)}+\parallel \nabla D(P^n)\parallel_{L^p(0,t;L^q)},\notag\\DE_n(t)=&\parallel D(d^n)\parallel_{L^\infty(0,t;B_{q,p}^{2(1-\frac{1}{p})})}+\parallel D(d^n)\parallel_{L^p(0,t;W^{2,q})}\\&+\parallel \partial_t D(d^n)\parallel_{L^p(0,t;L^q)},\notag\\
DH_n(t)=&DF_n(t)+DE_n(t).
\end{align} A direct computation shows that the triplet $(D(u^n),D(d^n),D(P^n))$ satisfies
\begin{equation}
 \left\{\begin{array}{ll}
\frac{\partial D(u^n)}{\partial t}-\triangle D(u^n)+\nabla D(P^n)\\=-u^n\cdot\nabla u^n+u^{n-1}\cdot \nabla u^{n-1}-\nabla(\frac{\mid \nabla d^n\mid^2}{2})\\+\nabla(\frac{\mid \nabla d^{n-1}\mid^2}{2})
-\triangle d^n\cdot\nabla d^n+\triangle d^{n-1}\cdot \nabla d^{n-1},
\ \ &in\ \  \mathbb{R}^{+}\times \Omega,\\
\\
\frac{\partial D(d^n)}{\partial t}-\triangle D(d^n)\\=-u^n\cdot \nabla d^n+u^{n-1}\cdot\nabla d^{n-1}\\
+\mid\nabla d^n\mid^2d^n-\mid\nabla d^{n-1}\mid^2d^{n-1}, \ \ &in\ \ \mathbb{R}^{+}\times \Omega,\\
\\
\nabla\cdot D(u^n)=0,\ \ \ \ \ \ \ \ \ \ \ \ \ \ \ \ \ \ \ \ \ \ \ &in\ \ \mathbb{R}^{+}\times \Omega, \\
\\
\int_{\Omega} D(P^n)dx=0,
\end{array}\right.
\end{equation}  with initial-boundary value
\begin{align}
D(u^n)|_{t=0}=D(u^n)|_{\partial \Omega}=0, \\
D(d^n)|_{t=0}=D(d^n)|_{\partial \Omega}=0.
\end{align} Using (2.12)-(2.13), we get
\begin{align}
&\parallel-u^n\cdot\nabla u^n+u^{n-1}\cdot\nabla u^{n-1}\parallel_{L^p(0,t;L^q)}\\\leq &\parallel u^n\cdot \nabla(D(u^{n-1}))\parallel_{L^p(0,t;L^q)}+\parallel \nabla u^{n-1}\cdot D(u^{n-1})\parallel_{L^p(0,t;L^q)}\notag
\\\leq &t^\frac{1}{p}\parallel u^n\parallel_{L^\infty(0,t;D_{A_q}^{1-\frac{1}{p},p})}\parallel D(u^{n-1})\parallel_{L^\infty(0,t;D_{A_q}^{1-\frac{1}{p},p})}\notag\\&+t^\frac{1}{p}\parallel u^{n-1}\parallel_{L^\infty(0,t;D_{A_q}^{1-\frac{1}{p},p})}\parallel D(u^{n-1})\parallel_{L^\infty(0,t;D_{A_q}^{1-\frac{1}{p},p})},\notag\\
&\parallel\triangle d^n\cdot\nabla d^n-\triangle d^{n-1}\cdot\nabla d^{n-1}\parallel_{L^p(0,t;L^q)}\\\leq
&\parallel  d^n\parallel_{L^\infty(0,t;B_{q,p}^{2(1-\frac{1}{p})})}\parallel D(d^{n-1})\parallel_{L^p(0,t;W^{2,q})}\notag\\&+\parallel  d^{n-1}\parallel_{L^p(0,t;W^{2,q})}\parallel D(d^{n-1})\parallel_{L^\infty(0,t;B_{q,p}^{2(1-\frac{1}{p})})},\notag\\
&\parallel \nabla^2 d^n\cdot \nabla d^n-\nabla^2 d^{n-1}\cdot \nabla d^{n-1}\parallel_{L^p(0,t;L^q)}\\\leq
&\parallel  d^n\parallel_{L^\infty(0,t;B_{q,p}^{2(1-\frac{1}{p})})}\parallel D(d^{n-1})\parallel_{L^p(0,t;W^{2,q})}\notag\\&+\parallel  d^{n-1}\parallel_{L^p(0,t;W^{2,q})}\parallel D(d^{n-1})\parallel_{L^\infty(0,t;B_{q,p}^{2(1-\frac{1}{p})})},\notag\\
&\parallel u^n\cdot \nabla d^n-u^{n-1}\cdot \nabla d^{n-1}\parallel_{L^p(0,t;L^q)}\\\leq &t^\frac{1}{p}\parallel u^n\parallel_{L^\infty(0,t;D_{A_q}^{1-\frac{1}{p},p})}\parallel D(d^{n-1})\parallel_{L^\infty(0,t;B_{q,p}^{2(1-\frac{1}{p})})}\notag\\&+t^\frac{1}{p}\parallel d^{n-1}\parallel_{L^\infty(0,t;B_{q,p}^{2(1-\frac{1}{p})})}\parallel D(u^{n-1})\parallel_{L^\infty(0,t;D_{A_q}^{1-\frac{1}{p},p})},\notag\\
&\parallel \mid\nabla d^n\mid^2d^n-\mid\nabla d^{n-1}\mid^2d^{n-1}\parallel_{L^p(0,t;L^q)}\\\leq&2t^\frac{1}{p}\parallel d^n\parallel^2_{L^\infty(0,t;B_{q,p}^{2(1-\frac{1}{p})})}\parallel D(d^{n-1})\parallel_{L^\infty(0,t;B_{q,p}^{2(1-\frac{1}{p})})}\notag\\&+t^\frac{1}{p}\parallel d^{n-1}\parallel_{L^\infty(0,t;B_{q,p}^{2(1-\frac{1}{p})})}\parallel D(d^{n-1})\parallel_{L^\infty(0,t;B_{q,p}^{B_{q,p}^{2(1-\frac{1}{p})}})}.\notag
\end{align} Using (3.21) and (3.29)-(3.33), we obtain that
\begin{align}
&DE_n(t)\leq Ct^\frac{1}{p}(DE_{n-1}(t)+DF_{n-1}(t)),\\
&DF_n(t)\leq C(t^\frac{1}{p}DF_{n-1}(t)+DE_{n-1}(t)).
\end{align} Thus, we get
\begin{align}
DH_{n+1}\leq 4 C^2(t^\frac{2}{p}+t^\frac{1}{p})H_{n-1}(t),
\end{align} where the constant $C$ only depends on $u_0,d_0$ and $\Omega$.
Choose $T_1\leq T_0$ such that $4C^2(T_1^\frac{2}{p}+T_1^\frac{1}{p})\leq \frac{1}{2}.$ Then
\begin{align}
DH_{n+1}(t)\leq\frac{1}{2} DH_{n-1}(t), \ \ \ on \ \ \ [0,T_1].
\end{align} It is clear that $\{(u^n,d^n,P^n)\}_{n\in\mathbb{N}}$ converges in $E^{q,p}_{T_1}.$
\end{proof}

Let $(u,d,p)$ be the limit of $\{(u^n,d^n,P^n)\}$, we prove that $(u,d,P)$ is a solution to the system
(3.1)-(3.2). It suffices to show that the items in the system (3.3) converge to corresponding items in $L^p(0,T_1;L^q).$ We only prove the convergence of $\nabla\cdot(\nabla d^n\odot\nabla d^n)$, the others can be proved similarly.
\begin{align}
&\parallel \nabla\cdot(\nabla d^n\odot\nabla d^n)-\nabla\cdot(\nabla d\odot\nabla d)\parallel_{L^p(0,T_1;L^q)}\\\leq &\parallel \nabla \frac{\mid \nabla d^n\mid^2}{2}-\nabla \frac{\mid \nabla d\mid^2}{2}\parallel_{L^p(0,T_1;L^q)}\notag\\&+
\parallel \triangle d^n\cdot \nabla d^n-\triangle d\cdot \nabla d\parallel_{L^p(0,T_1,L^q)}\notag\\\leq &2(\parallel d^n-d\parallel_{L^p(0,T_1;W^{2,q})}\parallel d^n\parallel_{L^\infty(0,T_1;B_{q,p}^{2(1-\frac{1}{p})})}\notag\\&+\parallel d\parallel_{L^p(0,T_1;W^{2,q})}\parallel d^n-d\parallel_{L^\infty(0,T_1;B_{q,p}^{2(1-\frac{1}{p})})}).\notag
\end{align} Using the convergence of $(u^n,d^n)$ and (3.21) we obtain
\begin{align}
&\parallel \nabla\cdot(\nabla d^n\odot\nabla d^n)-\nabla\cdot(\nabla d\odot\nabla d)\parallel_{L^p(0,T_1;L^q)}\rightarrow 0,\\
&\parallel u^n\cdot\nabla u^n-u\nabla u\parallel_{L^p(0,T_1;L^q)}\rightarrow 0,\\
&\parallel u^n\cdot\nabla d^n-u\cdot\nabla d\parallel_{L^p(0,T_1;L^q)}\rightarrow 0,\\
&\parallel \mid\nabla d^n \mid^2d^n-\mid\nabla d \mid^2d\parallel_{L^p(0,T_1;L^q)}\rightarrow 0.
\end{align} It is easy to see $\nabla u=0$ and $\int_\Omega Pdx=0$. As a direct consequence of Theorem 2.2, we get $\mid d\mid=1$. Thus, we have proved that the function triplet $(u,d,P)$ is a solution to the system (3.1)-(3.2).

Now, we prove the uniqueness. Let $(u^1,d^1,P^1)$ and $(u^2,d^2,P^2)$ be two solutions of the system (3.1) with the initial-boundary conditions (3.2). Denote
\begin{align*}
\delta u=u^1-u^2,\ \ \delta d=d^1-d^2,\ \ \delta P=P^1-P^2.
\end{align*} Then, the triplet $(\delta u,\delta d,\delta P)$ satisfies the following system
\begin{equation}
 \left\{\begin{array}{ll}
\frac{\partial \delta u}{\partial t}-\triangle \delta u+\nabla \delta P\\=-u^2\cdot\nabla u^2+u^1\cdot \nabla u^1-\nabla(\frac{\mid \nabla d^2\mid^2}{2})\\+\nabla(\frac{\mid \nabla d^1\mid^2}{2})
-\triangle d^2\cdot\nabla d^2+\triangle d^1\cdot \triangle d^1,
\ \ &in\ \  (0,T_1)\times \Omega,\\
\\
\frac{\partial \delta d}{\partial t}-\triangle \delta d\\=-u^2\cdot \nabla d^2+u^1\nabla d^1\\
+\mid\nabla d^2\mid^2d^2-\mid \nabla d^1\mid^2d^1, \ \ &in\ \ (0,T_1)\times \Omega,\\
\\
\nabla\cdot \delta u=0,\ \ \ \ \ \ \ \ \ \ \ \ \ \ \ \ \ \ \ \ \ \ \ &in\ \ (0,T_1)\times \Omega, \\
\\
\int_{\Omega} \delta Pdx=0
\end{array}\right.
\end{equation}  with initial-boundary value
\begin{align}
\delta u|_{t=0}=\delta u|_{\partial \Omega}=0, \\
\delta d|_{t=0}=\delta d|_{\partial \Omega}=0.
\end{align} Let
\begin{align}
G(t)=&\parallel \delta_n u\parallel_{L^\infty(0,t;D_{A_q}^{1-\frac{1}{p},p})}+\parallel \delta_n u\parallel_{L^p(0,t;W^{2,q})}+\parallel \partial_t \delta_n u\parallel_{L^p(0,t;L^q)}\\&+\parallel \delta_n d\parallel_{L^\infty(0,t;B_{q,p}^{2(1-\frac{1}{p})})}+\parallel \delta_n d\parallel_{L^p(0,t;W^{2,q})}+\parallel\partial_t\delta_n d\parallel_{L^p(0,t;L^q)}\notag\\&+\parallel \nabla \delta_n P\parallel_{L^p(0,t;L^q)}.\notag
\end{align} Repeating the argument from (3.22)-(3.37), we get
\begin{align}
G(t)\leq \frac{1}{2}G(t) \ \ \ on\ \ \ [0,T_1].
\end{align} Hence, $G(t)=0$ for $t\in[0,T_1]$, which implies the uniqueness on the interval $[0,T_1]$.

To complete this section, we show the continuous dependence. Noticing the proof of Lemma 3.1 and Lemma 3.2, more precisely, from (3.19) and (3.36), we can deduce that for any given initial-boundary data, there exists a $T_1>0$ only depending on the initial-boundary data such that $H(t)\leq2C(1+C^2H_0) H_0$ on $[0,T_1]$. Let
\begin{align} &(u_0,d_0)\in D_{A^q}^{1-\frac{1}{p},p}\times (B_{q,p}^{2(1-\frac{1}{p})}\cap C^{2,\alpha}(\partial \Omega)),\\ &(u^n_0,d^n_0)\rightarrow (u_0,d_0)\ \ in\ \ D_{A^q}^{1-\frac{1}{p},p}\times (B_{q,p}^{2(1-\frac{1}{p})}\cap C^{2,\alpha}(\partial \Omega)).\end{align} Assume that $(u^n,d^n)$ is the corresponding solution with the initial-boundary condition $(u^n_0,d^n_0)$ and  $(u,d)$ is the solution with initial-boundary data $(u_0,d_0)$. Define
\begin{align}
\delta_n u=u-u^n,\ \ \delta_n d=d-d^n,\ \ \delta_n P=P-P^n.
\end{align} Then, the triplet $(\delta_nu,\delta_nd,\delta_n P)$ satisfies the following system
\begin{equation}
 \left\{\begin{array}{ll}
\frac{\partial \delta_n u}{\partial t}-\triangle \delta_n u+\nabla \delta_n P\\=-u\cdot\nabla u+u^n\cdot \nabla u^n-\nabla(\frac{\mid \nabla d\mid^2}{2})\\+\nabla(\frac{\mid \nabla d^n\mid^2}{2})
-\triangle d\cdot\nabla d+\triangle d^n\cdot \triangle d^n,
\ \ &in\ \  (0,T_1)\times \Omega,\\
\\
\frac{\partial \delta_n d}{\partial t}-\triangle \delta_n d\\=-u\cdot \nabla d+u^n\nabla d^n\\
+\mid\nabla d\mid^2d-\mid \nabla d^n\mid^2d^n, \ \ &in\ \ (0,T_1)\times \Omega,\\
\\
\nabla\cdot \delta_n u=0,\ \ \ \ \ \ \ \ \ \ \ \ \ \ \ \ \ \ \ \ \ \ \ &in\ \ (0,T_1)\times \Omega, \\
\\
\int_{\Omega} \delta_n Pdx=0
\end{array}\right.
\end{equation} with initial-boundary conditions
\begin{align}
\delta_nu|_{\partial p Q_{T_1}}=u_0-u^n_0,\\
\delta_nd|_{\partial p Q_{T_1}}=d_0-d^n_0.
\end{align} Let
\begin{align}
G_n(t)=&\parallel \delta u\parallel_{L^\infty(0,t;D_{A_q}^{1-\frac{1}{p},p})}+\parallel \delta u\parallel_{L^p(0,t;W^{2,q})}+\parallel \partial_t \delta u\parallel_{L^p(0,t;L^q)}\\&+\parallel \delta d\parallel_{L^\infty(0,t;B_{q,p}^{2(1-\frac{1}{p})})}+\parallel \delta d\parallel_{L^p(0,t;W^{2,q})}+\parallel\partial_t\delta d\parallel_{L^p(0,t;L^q)}\notag\\&+\parallel \nabla \delta P\parallel_{L^p(0,t;L^q)}.\notag\\
G_n(0)=&\parallel u_0-u^n_0\parallel_{D_{A^q}^{1-\frac{1}{p},p}}+\parallel d_0-d^n_0\parallel_{B_{q,p}^{2(1-\frac{1}{p})}\cap C^{2,\alpha}(\partial \Omega)}.
\end{align} Using Lemma 2.2 and Theorem 2.1, and repeating the proof from (3.29) to (3.36), we have that
\begin{align}
G_n(t)\leq CG_n(0)+\frac{1}{2}G_n(t),\ \ on \ \ [0,T_1],
\end{align} where the constant $C$ only depends on $u_0$, $d_0$ and $T_1$. This implies $G_n(t)\leq C G_n(0)$.
Thus, we get the continuous dependence of the initial-boundary data.

Now, we complete the proof of well-posedness.

\section{Global existence }
\par
In this section, our aim is to extend the local solution which is established in Section 3 to a global solution in case of small initial-boundary conditions.

Given any unit vector $e\in S^2$ such that $d_0|\partial \Omega=e$. Let $(u,d)$ is a solution of the system
(3.1)-(3.2). A directly computation shows that $(u, d-e)$ satisfies
\begin{equation}
 \left\{\begin{array}{ll}
\frac{\partial u}{\partial t}+u\cdot\nabla u-\triangle u+\nabla P=-\nabla\cdot(\nabla d\bigodot\nabla d),
\ \ &in\ \  \mathbb{R}^{+}\times \Omega,\\
\frac{\partial d}{\partial t}+u\cdot\nabla d=(\triangle d+\mid\nabla d\mid^{2}d+e\mid\nabla d\mid^{2}), \ \ &in\ \ \mathbb{R}^{+}\times \Omega,\\
\nabla\cdot u=0,\ \ \ \ \ \ \ \ \ \ \ \ \ \ \ \ \ \ \ \ \ \ \ &in\ \ \mathbb{R}^{+}\times \Omega, \\
\int_{\Omega} Pdx=0,
\end{array}\right.
\end{equation} with initial-boundary conditions
\begin{align}
\left\{\begin{array}{ll}
(u(0,x),d(0,x))=(u_0(x),d_0(x)-e)\ \ x\in\Omega,\\
(u(t,x),d(t,x))=(0,0),\ \ \ \ \ (t,x)\in\mathbb{R}^{+}\times\partial\Omega.
\end{array}\right.
\end{align}

We first let $T^*< \infty$ be the maximal time of existence for $(u,d-e,P).$ Define
\begin{align}
H(t)=&\parallel u\parallel_{L^\infty(0,t;D_{A_q}^{1-\frac{1}{p},p})}+\parallel u\parallel_{L^p(0,t;W^{2,q})}+\parallel \partial_t u\parallel_{L^p(0,t;L^q)}\\&+\parallel d-e\parallel_{L^p(0,t;B_{q,p}^{2(1-\frac{1}{p})})}+\parallel \partial_t (d-e)\parallel_{L^p(0,t;L^q)}+\parallel P\parallel_{L^p(0,t;W^{1,q})}\notag\\&+\parallel d-e\parallel_{L^p(0,t;W^{2,q})},\notag\\
H_0=&\parallel u_0\parallel_{D_{A_q}^{1-\frac{1}{p},p}}+\parallel d_0-e\parallel_{B_{q,p}^{2(1-\frac{1}{p})}}.
\end{align} Using the system (3.1), Theorem 2.1 and Lemma 2.2, we get
\begin{align}
H(t)\leq &C(H_0+\parallel u\cdot\nabla u\parallel_{L^p(0,t;L^q)}+\parallel \nabla\cdot (\nabla (d-e)\odot\nabla (d-e))\parallel_{L^p(0,t;L^q)}\\&+\parallel u\cdot\nabla (d-e)\parallel_{L^p(0,t;L^q)}+\parallel \mid\nabla (d-e)\mid^2(d-e)\parallel_{L^p(0,t;L^q)}+\parallel\mid\nabla (d-e)\mid^2e\parallel_{L^p(0,t;L^q)}).\notag
\end{align} Noticing that $\mid d\mid=\mid e\mid=1$ and (2.12)-(2.13), we have that
\begin{align}
H(t)\leq &C(H_0+\parallel u\parallel_{L^\infty(0,t;D_{A_q}^{2(1-\frac{1}{p})})}\parallel u\parallel_{L^p(0,t;W^{2,q})}\\&+\parallel d-e\parallel_{L^\infty(0,t;B_{q,p}^{2(1-\frac{1}{p})})}\parallel d-e\parallel_{L^p(0,t;W^{2,q})}\notag\\&+\parallel u\parallel_{L^\infty(0,t;D_{A_q}^{1-\frac{1}{p},p})}\parallel d-e\parallel_{L^p(0,t;W^{2,q})}\notag\\&+\parallel d-e\parallel_{L^\infty(0,t;B_{q,p}^{2(1-\frac{1}{p})})}\parallel d-e\parallel_{L^p(0,t;W^{2,q})}).\notag
\end{align} From the definition of $H(t)$ and (4.4), we have
\begin{align}
H(t)\leq C(H_0+4H^2(t)), \ \ on \ \ [0,T^*).
\end{align} If we take $H_0\leq\frac{1}{32C^2}$, then the continuity and monotonicity of $H(t)$ yields that there exists a positive
$T$ such that
\begin{align}
H(t)\leq 2CH_0\leq \frac{1}{16C}, \ \ on \ \ [0,T].
\end{align}  On the other hand, (4.5) yields that
\begin{align}
H(t)\leq \frac{1-(1-16C^2H_0)^{\frac{1}{2}}}{8C}\ \ or\ \ H(t)\geq \frac{1+(1-16C^2H_0)^{\frac{1}{2}}}{8C}\ \ on\ \ [0,T^*).
\end{align}
Since $H(t)\leq \frac{1}{16C^2}$ on $[0,T]$, the continuity of $H(t)$ implies that
\begin{align}
H(t)\leq \frac{1-(1-16C^2H_0)^{\frac{1}{2}}}{8C}.
\end{align} This yields that
\begin{align}
\parallel(u,d-e,P)\parallel_{E^{p,q}_{T^*}}\leq\frac{1}{8C}.
\end{align} It is impossible and the solution exists globally in time. This implies that $(u,d,P)$ is a global solution of the system (3.1)-(3.2).

The proof of Theorem 1.1. is complete.

\bigskip
\noindent\textbf{Acknowledgments} This work was partially
supported by NNSFC (No. 10971235), RFDP (No. 200805580014),
NCET-08-0579 and the key project of Sun Yat-sen University.

\end{document}